\documentclass[11pt]{article}
\usepackage{amssymb,amsmath,latexsym,graphicx}
\newenvironment{proof}{\textit{Proof.}}{\hfill $\Box$}

\newcommand{\comment}[1]{}
\newfont{\menutt}{cmtt8}

\newcommand{\cross}{\mathrm{cross}}
\newcommand{\nest}{\mathrm{nest}}
\newcommand{\leftright}{left-right }

\newtheorem{thm}{Theorem}[section]

\newtheorem{cor}[thm]{Corollary}
\newtheorem{prop}[thm]{Proposition}
\newtheorem{lemma}[thm]{Lemma}

\title{$k$-noncrossing and $k$-nonnesting graphs and fillings of Ferrers diagrams}
\author{Anna de Mier\thanks{Department of Applied Mathematics and Institute of
Theoretical Computer Science, Charles
University, Malostransk\'e n\'am. 25, 118 00 Praha 1, Czech Republic. \emph{email address:} \texttt{demier@kam.mff.cuni.cz}}}
\date{}

\begin{document}

\maketitle

\begin{abstract}
 We give a correspondence between graphs with a given degree
sequence and fillings of Ferrers diagrams by nonnegative integers
with prescribed row and column sums. In this setting,
$k$-crossings and $k$-nestings of the graph become occurrences of
the identity and the antiidentity matrices in the filling. We use
this to show the equality of the numbers of $k$-noncrossing and
$k$-nonnesting graphs with a given degree sequence. This
generalizes the analogous result for matchings and partition
graphs of Chen, Deng, Du, Stanley, and Yan, and extends results of
Klazar to
 $k>2$. Moreover,
this correspondence reinforces the links recently discovered by
Krattenthaler between fillings of diagrams and the results of Chen
et al.

\end{abstract}


\section{Introduction}\label{sec:intro}


Let $G$ be a graph on $[n]$; unless otherwise stated, we allow
multiple edges and isolated vertices, but no loops. Two edges
$\{i,j\}$ and $\{k,l\}$ are a \emph{crossing} if $i<k<j<l$ and
they are a \emph{nesting} if $i<k<l<j$. If we draw the vertices of
$G$ on a line and represent the corresponding edges by arcs above
the line, crossings and nestings have the obvious geometric
meaning. A graph without crossings (respectively, nestings) is
called \emph{noncrossing} (resp., \emph{nonnesting}).
Klazar~\cite{klazar_crossings} proves the equality between the
numbers of noncrossing and nonnesting simple graphs, counted by
order, and also between the numbers of noncrossing and nonnesting
graphs without isolated vertices, counted by size. The purpose of
this paper is to study analogous results for sets of $k$ pairwise
crossing and $k$ pairwise nested edges.

A \emph{$k$-crossing} is a set of $k$ edges every two of them
being a crossing, that is, edges $\{i_1,j_1\},\ldots,$
$\{i_k,j_k\}$ such that $i_1<i_2<\cdots<i_k<j_1<\cdots< j_k$. A
 \emph{$k$-nesting} is a set of $k$ edges pairwise
nested, that is, $\{i_1,j_1\},\ldots,$ $\{i_k,j_k\}$ such that
$i_1<i_2<\cdots<i_k<j_k<\cdots< j_1$.
 A graph with no $k$-crossing is called
\emph{$k$-noncrossing} and  a graph with no $k$-nesting is called
\emph{$k$-nonnesting}. The largest $k$ for which a graph $G$ has a
$k$-crossing (respectively, a $k$-nesting) is denoted $\cross(G)$
(resp., $\nest(G)$). The aim of this paper is to show that the
number of $k$-noncrossing graphs equals the number of
$k$-nonnesting graphs, counted by  order, size, and degree
sequences. This problem was originally posed by Martin Klazar and
we learned of it at the Homonolo 2005 workshop~\cite{homonolo};
the case where the number of vertices of the graph is $2k+1$ was
proved by A. P\'or (unpublished). Our main result
(Theorem~\ref{thm:kcrossknest}) states that the numbers of
$k$-noncrossing and $k$-nonnesting graphs with a given degree
sequence are the same.

Chen~et~al.~\cite{match_and_part} prove the equality  of the
numbers  of $k$-noncrossing and $k$-nonnesting graphs  for two
subclasses of graphs, namely for perfect matchings and for
partition graphs, also counted by degree sequences (under a
different but equivalent terminology). A perfect matching is a
graph where each vertex has degree one, and a partition graph is a
graph that is a disjoint union of monotone paths, that is, where
each vertex has at most one edge to its right and at most one to
its left. The latter correspond in a natural way to set
partitions, hence the result can be stated in terms of these. The
paper~\cite{match_and_part} also contains other identities and
enumerative results on $k$-noncrossing and $k$-nonnesting
matchings and partitions. Krattenthaler~\cite{krattenthaler}
deduces most of these from his more general results on fillings of
Ferrers diagrams. In this paper we also use fillings of diagrams
to prove results about graphs. The difference is that whereas
in~\cite{krattenthaler}, and also in~\cite{jonsson}, the results
about graphs follow from general theorems by restricting the shape
of the diagram, here we show that the results about graphs are in
fact equivalent to those about fillings with arbitrary shapes.

The main idea is to encode graphs by fillings of Ferrers diagrams
in such a way that $k$-crossings and $k$-nestings are easy to
recognize. A $k$-noncrossing ($k$-nonnesting)
 graph becomes a filling of a diagram
that avoids the identity (antiidentity) matrix of order $k$, and
the degree sequence of the graph can be recovered from the shape
of the diagram and the row and column sums of the filling. Then
proving that there are as many $k$-noncrossing as $k$-nonnesting
graphs is equivalent to showing that the numbers of fillings
avoiding these two matrices are the same.  This idea generalizes
easily to other subgraphs in addition to crossing and nestings,
and allows us to show that the study of fillings of Ferrers
diagrams with forbidden configurations is equivalent to the study
of graphs avoiding certain subgraphs, in the sense defined in
Section~\ref{sec:degree}.

The structure of the paper is as follows. In
Section~\ref{sec:general} we show that the equality of the numbers
of $k$-noncrossing and $k$-nonnesting graphs counted by size and
order is already in the literature, although not explicitly stated
in this form. We introduce some notation on pattern avoiding
fillings of Ferrers diagrams and we rephrase results of
Krattenthaler~\cite{krattenthaler} and Jonsson and
Welker~\cite{pfaffians}
 in terms of  $k$-noncrossing and
$k$-nonnesting graphs. Section~\ref{sec:degree} introduces a new
correspondence between graphs and fillings of diagrams that keeps
track of degree sequences.  Then we discuss why, from the
perspective of pattern avoiding, graphs and fillings of diagrams
are equivalent objects. In particular, showing that the number of
$k$-noncrossing graphs with a fixed degree sequence equals the
number of such $k$-nonnesting graphs is equivalent to proving a
result on fillings of diagrams with restrictions on the row and
column sums.
 Our proof is an adaptation of the one
in~\cite{filling_boards} to allow arbitrary entries in the
filling, and this is the content of Section~\ref{sec:proof}.  We
conclude with some remarks and open questions.


\section{Fillings of diagrams}\label{sec:general}

\begin{figure}[ht]
\begin{center}
\includegraphics[width=9cm]{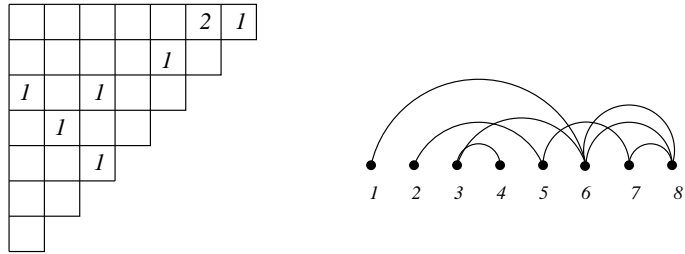}
\caption{Left: a filling of a diagram that contains $I_3$ and
$J_2$ but avoids $J_3$. Right: the graph determined by the filling
has a $3$-nesting and several $2$-crossings, but no
$3$-crossing.}\label{fig:diagram1}
\end{center}
\end{figure}

We start by setting some notation on fillings of Ferrers diagrams.
Let $\lambda=(\lambda_1,\lambda_2,\ldots,\lambda_k)$ be an integer
partition. The \emph{Ferrers diagram of shape $\lambda$} (or
simply a \emph{diagram}) is the arrangement of square cells,
left-justified and from top to bottom, having $\lambda_i$ cells in
row $i$, for $i$ with $1\leq i \leq k$. For a Ferrers diagram $T$
of shape $\lambda$
 with rows indexed from top to bottom and columns
from left to right,  a \emph{filling}  $L$ of $T$ consists of
assigning a nonnegative integer to each cell of the diagram. We
say that a cell is \emph{empty} if it has been assigned the
integer $0$.  Let $M$ be an $s\times t$ $0-1$ matrix. We say that
the filling \emph{contains $M$} if there is a selection of rows
$(r_1,\ldots,r_s)$ and columns $(c_1,\ldots,c_t)$ of $T$ such that
if $M_{i,j}=1$ then the cell $(r_i,c_j)$ of $T$ is nonempty and
moreover the cell $(r_s,c_t)$ is in the diagram (in other words,
we require that the matrix $M$ is fully contained in $T$). We say
that the filling \emph{avoids $M$} if there is no such selection
of rows and columns. If a filling $L$ contains $M$, by an
\emph{occurrence} of $M$ we mean the set of cells of $T$ that
correspond to the $1$'s in $M$.
 We are mainly concerned
about diagrams avoiding the identity matrix $I_t$ and the
antiidentity matrix $J_t$; the latter is the matrix with $1$'s in
the main antidiagonal and $0$'s elsewhere.  As an example of these
concepts, Figure~\ref{fig:diagram1} shows a filling of a diagram
of shape $(7,6,5,4,3,2,1)$ that contains the matrices $I_3$ and
$J_2$ but avoids $J_3$. (For clarity, we omit the zeros
corresponding to the empty cells.)

Studying fillings of diagrams avoiding matrices  is a natural
generalization of pattern avoiding permutations, as explained
in~\cite{filling_boards, stankova_west}. We  explore two types of
connections between graphs and fillings of diagrams. The first one
is straightforward, being essentially the adjacency matrix, and it
has been used in~\cite{jonsson, krattenthaler} to derive results
on $k$-noncrossing maximal graphs and $k$-noncrossing and
$k$-nonnesting matchings and partitions.

Suppose $G$ is a graph on $[n]$ and consider a diagram $\Delta$ of
shape $(n-1,n-2,\ldots,2,1)$. Then if there are $d\geq 0$ edges
joining vertices $i$ and $j$, with $i<j$, fill the cell of column
$i$ and row $n-j+1$ with $d$. Let this filling of the diagram be
called $\Delta(G)$. Obviously the sum of the entries of
$\Delta(G)$ is the number of edges of $G$ and the number of
vertices is just one plus the number of rows of $\Delta$. If $G$
is a simple graph, then $\Delta(G)$ is a $0-1$ filling. If the
edges $\{i_1,j_1\},\ldots,\{i_k,j_k\}$ are a $k$-nesting of $G$,
then
 $\Delta(G)$ contains the $k\times k$ identity matrix $I_k$ in columns
$i_1,\ldots,i_k$ and rows $n-j_1+1,\ldots,n-j_k+1$. Similarly, if
$G$ contains a $k$-crossing, then $\Delta(G)$ contains the
antiidentity matrix $J_k$ (the condition $i_k<j_1$ guarantees that
the matrix is indeed contained in the diagram).

Krattenthaler~\cite{krattenthaler} derives many of the results of
Chen et al.~\cite{match_and_part} for matchings and partitions by
specializing to $\Delta(G)$
 his results on fillings of diagrams avoiding large identity
or antiidentity matrices. His Theorem~13 gives a generalization to
arbitrary graphs which is implicitly included in the remark after
it; we explicitly state his result here (see also the comment
after Theorem~\ref{thm:kcrossknest}
 in the next section). The following
is a weaker version of~\cite[Theorem 13]{krattenthaler}

\begin{thm}\label{thm:thm13}
For any diagram $T$ and any integer $m$, consider
 fillings of $T$ with nonnegative integers adding up to $m$. Then
for each $k>1$, the number
of such fillings that do not contain the identity matrix $I_k$ equals the
number of fillings that do not contain the antiidentity matrix $J_k$.
\end{thm}

By restricting to $T=\Delta$ we immediately get the following.

\begin{cor}\label{cor:kcrossknestmulti}
The number of $k$-noncrossing graphs with $n$ vertices and $m$ edges equals the
number of $k$-nonnesting such graphs.
\end{cor}

Actually, from the statement of~\cite[Theorem 13]{krattenthaler}
one gets a stronger result. For this we need to introduce weak
$k$-crossings and weak $k$-nestings. The edges
$\{i_1,j_1\},\ldots,\{i_k,j_k\}$ are a \emph{weak $k$-crossing} if
$i_1\leq i_2\leq \cdots \leq i_k < j_1\leq \cdots\leq j_k$;
similarly, they are a \emph{weak $k$-nesting} if $i_1\leq i_2\leq
\cdots \leq i_k <$ $j_k\leq \cdots\leq j_1$. Let $\cross^*(G)$
(respectively, $\nest^*(G)$) be the largest $k$ for which $G$ has
a weak $k$-crossing (resp., weak $k$-nesting). Then the following
is a corollary of the full version of~\cite[Theorem
13]{krattenthaler}.

\begin{cor}\label{cor:weak}
The number of graphs with $n$ vertices and $m$ edges with $\cross(G)=r$
and $\nest^*(G)=s$ equals the number of such graphs with $\cross^*(G)=s$
and $\nest(G)=r$.
\end{cor}

Ideally, one would like to have an analogous result proving the
symmetry of the distribution of $\cross(G)$ and $\nest(G)$ for all
graphs. This is known to be true for matchings and partition
graphs~\cite[Theorem 1 and Corollary 4]{match_and_part}; actually,
for these graphs weak crossings (respectively, weak nestings) are
the same as crossings (resp., nestings). For simple graphs, the
result would follow if Problem~2 in~\cite{krattenthaler} has a
positive answer for the diagram $\Delta$.

The bijection used to prove~\cite[Theorem 13]{krattenthaler} does
not preserve the values of the entries of the filling, so we
cannot deduce from it the corresponding result for simple graphs.
However, this follows from a  result of Jonsson and Welker. They
deal with fillings not of diagrams, but of stack polyominoes. A
\emph{stack polyomino} consists of taking a diagram, reflecting it
through the vertical axis, and gluing it to another (unreflected)
diagram. The content of a stack polyomino is the multiset of the
lengths of its columns. The definitions of fillings and
containment of matrices in stack polyominoes are analogous to
those for diagrams. The following is Corollary~6.5
of~\cite{pfaffians}.
 (The particular
case where $m$ below is maximal was proved in~\cite{jonsson}.)

\begin{thm}\label{thm:pfaffians}
The number of $0-1$ fillings of a stack polyomino with $m$ nonzero
entries that avoid the matrix $I_k$ depends only on the content of
the  polyomino and not on the ordering of the columns.
\end{thm}

  By a simple reflection argument we
get the following for the triangular diagram $\Delta$ of shape
$(n-1,n-2,\ldots,2,1)$: the number of $0-1$ fillings of $\Delta$
with $m$ non-zero entries and that avoid the matrix $I_k$ is the
same as those that avoid the matrix $J_k$. Hence we have the
following in terms of graphs.

\begin{cor}\label{cor:simple}
The number of $k$-noncrossing simple graphs on $n$ vertices
and $m$ edges equals the number of such $k$-nonnesting simple graphs.
\end{cor}

In the next section we deal with graphs with a  fixed degree
sequence. For this we  need to consider diagrams of arbitrary shapes,
since the correspondence between graphs and fillings
 is no longer restricted to  the
triangular diagram $\Delta$.


\section{Degree sequences and fillings with prescribed row and column sums}\label{sec:degree}


The \emph{left-right degree sequence} of a graph on $[n]$ is the
sequence $\left((l_i  ,  r_i)\right)_{1\leq i \leq n}$, where
$l_i$ (resp., $r_i$) is the left (resp., right) degree of vertex
$i$; by the left (resp., right) degree of $i$ we mean  the number
of edges that join $i$ to a vertex $j$ with $j<i$ (resp., $j>i$).
Obviously $l_i+r_i$ is the degree of vertex $i$ (loops are not
allowed). For instance, if $r_i\leq 1$ and $l_i\leq 1$ for all
$i$, then
 the graph is either a matching or a partition graph, perhaps with
some isolated vertices. If a graph $G$ has $D$ as its left-right
degree sequence, we say that $G$ is a graph \emph{on $D$}. A
useful way of thinking of left-right degree sequences is drawing
for each vertex $i$, $l_i$ half-edges going left and $r_i$
half-edges going right. Then a graph is just a way of matching
these half-edges; recall that we allow multiple edges. For
completeness we mention here that a sequence $\left((l_i,
r_i)\right)_{1\leq i \leq n}$  is the left-right degree sequence
of some graph on $[n]$ if and only if
\begin{equation}\label{eq:deg} \sum_{i=1}^n l_i = \sum_{i=1}^{n}
r_i \mbox{ and } \sum_{i=1}^k l_i \leq \sum_{i=1}^{k-1} r_i, \quad
\forall k\in [n].\end{equation}

 This and
the next section are devoted to proving that for each left-right
degree sequence $D$ there are as many $k$-noncrossing graphs on
$D$ as $k$-nonnesting.  We stress that the fact that we allow
multiple edges is essential, since if we restrict to simple graphs
the result does not hold. For instance, one can check that there
is one simple nonnesting graph with left-right degree sequence
$(0,2),(0,2),(1,1),(2,0),(2,0)$, but no such noncrossing simple
graph. However, it turns out that there is a bijection between
$k$-noncrossing and $k$-nonnesting simple graphs that preserves
left degrees (or  right degrees, but not both simultaneously).
This follows from the following result of Rubey~\cite[Theorem
4.2]{rubey}  applied to the filling $\Delta(G)$ by noting that the
sum of the entries in row $n-j+1$ of $\Delta(G)$ corresponds to
the left degree of vertex $j$. Rubey's result is for moon
polyominoes, but we state the version for stack polyominoes. (A
weaker version of this result was proved by
Jonsson~\cite[Corollary 26]{jonsson}.)

\begin{thm}\label{thm:cor26}
For any stack polyomino $\Lambda$ with $s$ rows and for any
sequence $(d_1,\ldots, d_s)$ of nonnegative integers, the number
of $0-1$ fillings of $\Lambda$ that avoid $I_k$ and have $d_i$
nonzero entries in row $i$ depends only on the content of
$\Lambda$ and not on the ordering of the columns.
\end{thm}

By the same reflection argument as at the end of Section~\ref{sec:general}
we obtain the following corollary.

\begin{cor}\label{cor:simplemax}
Let  $(l_2,\ldots,l_n)$ be a sequence of nonnegative integers.
Then the number of $k$-noncrossing simple graphs on $[n]$ with
vertex $i$ having left degree $l_i$ for $2\leq i\leq n$ is the
same as the number of such $k$-nonnesting simple graphs.
\end{cor}

 The main result of this paper says that by allowing multiple
 edges we can simultaneously fix left and right degrees.

\begin{thm}\label{thm:kcrossknest}
For any left-right degree sequence $D$, the number of
$k$-noncrossing graphs on $D$ equals the number of $k$-nonnesting
 graphs on $D$.
\end{thm}

This result generalizes to arbitrary graphs some of the results
of~\cite{match_and_part}, which are only for partition graphs but
also taking into account degree sequences (with different
terminology). To approach Theorem~\ref{thm:kcrossknest}
 we could use again the filling $\Delta(G)$
of the previous section and fix the sums of the entries in each
row and column. In this setting Theorem~\ref{thm:kcrossknest} is
again implicitly included in the remark after~\cite[Theorem
13]{krattenthaler} by keeping track of the changes in the
partitions involved in the proof of that theorem. (I am grateful
to Christian Krattenthaler for this observation.) However, our
approach consists of encoding graphs not by the triangular diagram
$\Delta$ but by an arbitrary diagram whose shape depends on the
degree sequence. By doing this we actually show that not only
results on $k$-noncrossing and $k$-nonnesting graphs can be
deduced from results on fillings of Ferrers diagrams avoiding
$I_k$ and $J_k$, but that actually these two families of results
are completely equivalent. Moreover,  we have an analogous
assertion for arbitrary matrices (see
Theorem~\ref{thm:equivalence}).

We start with an easy lemma that follows immediately from the fact
that the edges in a $k$-crossing or a $k$-nesting must be
vertex-disjoint.

\begin{lemma}\label{lem:split}
The number of $k$-noncrossing (resp., $k$-nonnesting) graphs with
left-right degree sequence $$(l_1, r_1),\ldots,(l_n,r_n)$$
is the same as that of those with left-right degree sequence
$$(l_1,r_1),\ldots,(l_{i-1},r_{i-1}),(l_i,0),(0,r_i),(l_{i+1},r_{i+1}),
\ldots,(l_n,r_n), $$ for any $i$ with $1\leq i \leq n$.
\end{lemma}

Hence it is enough to prove Theorem~\ref{thm:kcrossknest} for
left-right degree sequences whose elements $(l_i,r_i)$ are such
that either $l_i$ or $r_i$ is $0$. We call these graphs
\emph{left-right} graphs; note though that we do not require that
the degrees of the vertices alternate between right and left. The
case where both left and right degrees are $0$ corresponds to an
isolated vertex. \comment{; even if isolated vertices are easy to
deal with, allowing them makes some proofs (things,
 statements)  more
homogeneous, hence we keep them.}

\begin{figure}[ht]
\begin{center}
\includegraphics[width=9cm]{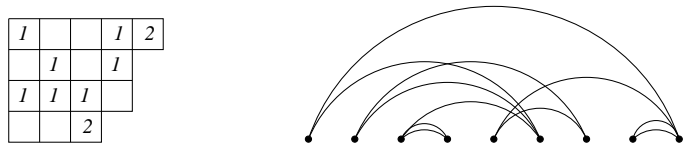}
\caption{A filling $L$ of a diagram with row sums $4,2,3,2$ and
column sums $2,2,3,2,2,$ and the corresponding graph $G(L)$.
}\label{fig:diagram2}
\end{center}
\end{figure}

We now describe a bijection between \leftright graphs
 and fillings of Ferrers diagrams of
arbitrary shape; this bijection has the property that the
left-right degree sequence of the graph can be recovered from the
shape and filling of the diagram. Let $G$ be a \leftright graph.
If the degree of vertex $i$ is of the form $(0,r_i)$ we say that
$i$ is \emph{opening}, and if it is of the form $(l_i,0)$ we say
that $i$ is \emph{closing}. An isolated vertex is both opening and
closing. Let $i_1,\ldots, i_c$ be the closing vertices of $G$ and
let $j_1,\ldots, j_o$ be the opening ones. For each closing vertex
$i$, let $p(i)$ be the number of vertices $j$ with $j<i$ that are
opening. We consider a diagram $T(G)$ of
 shape $(p(i_c),p(i_{c-1}),\ldots,p(i_1))$, and if there are $d$ edges
going from the opening vertex $j_s$ to the closing vertex $i_r$,
we fill
 the cell in column $s$ and row $c-r+1$ with the integer $d$
(see Figure~\ref{fig:diagram2}). Thus graphs with left degrees
$l_1,\ldots,l_c$ and right degrees $r_1,\ldots,r_o$ correspond to
fillings of this diagram with nonnegative entries such that the
sum of the entries in row $i$ is $l_i$ and the sum of the entries
in column $j$ is $r_j$. Conversely, any filling of a diagram
arises in this way. Indeed, given a filling $L$ of a diagram $T$,
the shape of $T$ gives the ordering of the opening and closing
vertices of the graph, the row and column sums give the left and
right degrees (it is easy to see that they must satisfy
equation~(\ref{eq:deg})), and the entries of the filling give the
edges of the graph.  Given a graph $G$, we denote by $L(G)$ the
filling of $T(G)$ corresponding to $G$. Similarly, given a filling
$L$ of a diagram, we denote by $G(L)$ the \leftright graph
corresponding to this filling.

In this
setting, it is immediate to check that again
 $k$-crossings of $G$ correspond to occurrences of $I_k$ in $L(G)$ and $k$-nestings to
occurrences of $J_k$.

By a \emph{diagram with prescribed row and column sums} we mean a
diagram and two sequences $(\rho_i)$ and $(\gamma_j)$  of
nonnegative integers such that the only fillings allowed for this
diagram are those where the row and column sums are given by the
sequences $(\rho_i)$ and $(\gamma_j)$. Given two matrices $M$ and
$N$, we say that they are \emph{equirestrictive} if for all
diagrams $T$ with prescribed row and column sums, the number of
fillings of $T$ that avoid $M$ equals the number of fillings of
$T$ that avoid $N$. With this notation,
Theorem~\ref{thm:kcrossknest} is an immediate consequence of the
following result, the proof of which is the content of the next
section.

\begin{thm}\label{thm:avoidk}
The  identity matrix $I_k$
and the antiidentity matrix $J_k$ are equirestrictive.
\end{thm}

\comment{ Let us mention that recently M. Rubey~\cite{rubey} has
proved the analogous result for fillings of moon polyominoes, also
with prescribed row and column sums. ASK MARTIN IF HE ALSO HAS
COLUMN SUMS FIXED, IT IS NOT CLEAR FROM HIS MESSAGE. }

Before moving to the proof of Theorem~\ref{thm:avoidk}, let us
make some remarks and point out some consequences of the proof.
We start by further exploring the bijection between \leftright
graphs and fillings
of diagrams.

Let $G$ and $H$ be graphs on $[n]$ and $[h]$, respectively, with
$h\leq n$. For the rest of this section we assume that $H$ is
simple (but $G$ can have multiple edges as usual). We say that $G$
contains $H$ if there is an order-preserving injection $\sigma:
[h]\rightarrow [n]$ such that if $\{i,j\}$ is an edge of $H$ then
$\{\sigma(i),\sigma(j)\}$ is an edge of $G$. For instance, a
$k$-noncrossing graph is a graph that does not contain the graph
on $[2k]$ with edges $\{1,k+1\},\{2,k+2\},\ldots,\{k,2k\}$.

A $0-1$ matrix $M$ with $s$ rows and $t$ columns can also be
viewed as a filling of the diagram of shape
$(t,t,\stackrel{(s)}{\ldots},t)$. By the correspondence between
graphs and fillings of diagrams described above, we have that $M$
gives a graph $G(M)$ with $t$ opening vertices and $s$ closing
vertices and such that all opening vertices appear before the
closing vertices. Let us call such a graph a \emph{split graph}, a
particular case being the graph of a $k$-crossing or a
$k$-nesting. As a consequence of the previous discussion we have
that in terms of containment of substructures (matrices or split
graphs), fillings of diagrams and graphs are equivalent objects.

\begin{thm}\label{thm:bijection}
For any split graph $H$ there is a  matrix $M(H)$ such that a
\leftright graph $G$ contains $H$ if and only if the filling
$L(G)$ contains $M(H)$. And conversely, for each matrix $M$ there
is a split graph $H(M)$ such that a filling $L$ of a diagram $T$
contains $M$ if and only if the graph $G(L)$ contains $H(M)$.
\end{thm}

Observe now that Lemma~\ref{lem:split} can be generalized by
substituting ``$k$-noncrossing graphs'' with ``graphs that do not
contain the split graph $H$''. Hence the following.

\begin{thm}\label{thm:equivalence}
Let $H$ and $H'$ be two split graphs. Then for any left-right
degree sequence $D$ there are as many graphs on
$D$ avoiding $H$ as graphs on $D$ avoiding $H'$ if and only if for each diagram
with prescribed row and column sums there are as many fillings avoiding
$M(H)$ as fillings avoiding $M(H')$.
\end{thm}

Following the notation for matrices, we say that two split graphs
$H$ and $H'$ are \emph{equirestrictive} if for any left-right
degree sequence $D$, there are as many graphs on $D$ avoiding $H$
as graphs on $D$ avoiding $H'$. All the split graphs that are
known to be equirestrictive are obtained from the graph of a
$k$-crossing or a $k$-nesting by using
Proposition~\ref{prop:avoidextension} from the next section. This
proposition states that if $M$ and $N$ are equirestrictive
matrices, then for any other matrix $A$ the matrices
$$\left( \begin{array}{cc} M & 0 \\ 0 & A \end{array}\right) \quad
\mbox{and} \quad  \left( \begin{array}{cc} N & 0 \\ 0 & A
\end{array}\right)$$ defined by blocks are also equirestrictive.
This has the following implications in terms of graphs. Given a
split graph $H$ on $[h]$, a \emph{$(k-H)$-crossing} is a graph on
$2k+h$ such that the graph induced by the vertices $[k]\cup
\{k+h+1,\ldots,2k+h\}$ is a $k$-crossing, the graph induced by
$\{k+1,\ldots,k+h\}$ is $H$, and there are no other edges. A
\emph{$(k-H)$-nesting} is defined similarly. Then by combining
Theorems~\ref{thm:avoidk} and~\ref{thm:equivalence} and
Proposition~\ref{prop:avoidextension} we deduce the following.

\begin{cor}\label{cor:kH}
For any split graph $H$ and any nonnegative integer $k$,
$(k-H)$-crossings and $(k-H)$-nestings are equirestrictive.
\end{cor}

Observe that if we take $H$ to be an $h$-nesting, a
$(k-H)$-nesting is a $(k+h)$-nesting, so it follows that a
$k$-nesting, a $k$-crossing, and any combination of a $t$-crossing
``over" a $(k-t)$-nesting are equirestrictive. However, it is not
true that  $t$-nestings over $(k-t)$-crossings are
equirestrictive, not even within matchings, as observed in the
remark after Theorem~1 of~\cite{jelinek_matchings}. This implies
also that there is no analogous version of
Proposition~\ref{prop:avoidextension} where $A$ is the top-left
block and $M$ and $N$ are the bottom-right blocks of the matrix.

Finally, we comment on the results known for $0-1$ fillings of
diagrams with row and column sums equal to $1$. Our correspondence
translates these results into results for matchings and partition
graphs, as we next explain.
 In the literature,
two permutation matrices $M$ and $N$ are called
\emph{shape-Wilf-equivalent} if for each diagram $T$ with row and
column sums set to $1$, the number of fillings avoiding $M$ equals
the number of fillings avoiding $N$. (In view of this notation, we
could have chosen the name graph-Wilf-equivalent instead of
equirestrictive.) Let $P$ be a $t\times t$ permutation matrix. The
split graph corresponding to $P$ is a matching (these are
sometimes called \emph{permutation matchings}). Now if two
permutation matrices $P$ and $P'$ are shape-Wilf-equivalent, then
by straightforward application of Theorem~\ref{thm:equivalence} we
have that for all graphs whose left and right degrees are one, the
number of graphs avoiding the matching $H(P)$ equals the number of
graphs avoiding the matching  $H(P')$. Since graphs with left and
right degrees one are exactly partition graphs, it turns out that
shape-Wilf-equivalence is equivalent to the matchings $H(P)$ and
$H(P')$ being equirestrictive among partition graphs, counted by
left-right degree sequences.

There are not many pairs of permutation matrices known to be
shape-Wilf-equivalent. Backelin, West, and
Xin~\cite{filling_boards} show that $I_k$ and $J_k$ are
shape-Wilf-equivalent; in graph theoretic terms, this gives
another alternative proof of the equality between $k$-noncrossing
and $k$-nonnesting partition graphs from~\cite{match_and_part}.
Let us mention here that Krattenthaler~\cite{krattenthaler}
deduces both the result of Chen et al. and that of Backelin, West,
and Xin from his Theorem~3, but for the first one he sets
$T=\Delta$ and for the second he restricts the number of non-empty
cells in the filling (and takes arbitrary shapes). Since these two
apparently unrelated results are in fact equivalent,  it is
obvious that they must follow from the same theorem,  but it is
interesting that they do in different ways. Another observation is
that by Lemma~\ref{lem:split}, and its generalization to split
graphs, if we know that two split graphs are equirestrictive
within matchings, then they are so within partition graphs. For
instance, a bijective proof of the equality of the numbers of
$k$-noncrossing and $k$-nonnesting matchings would immediately
give a bijection for $k$-noncrossing and $k$-nonnesting partition
graphs.

In addition to the matrices $I_t$ and $J_t$ and the ones that
follow from Proposition~\ref{prop:avoidextension}, the only other
pair of matrices known to be shape-Wilf-equivalent are
(see~\cite{stankova_west})
$$M(213)=\left( \begin{array}{ccc} 0& 0& 1 \\ 1&0&0 \\ 0&1&0 \end{array}
\right) \qquad \mathrm{and} \qquad M(132)=\left(\begin{array}{ccc}
0&1&0\\ 0&0&1\\1&0&0 \end{array}
 \right). $$
The graph theoretic version of this result has been independently
proved by Jel\'inek~\cite{jelinek_matchings}. It is not known to
us if $M(231)$ and $M(132)$ are also equirestrictive, or more
generally if there is a pair of shape-Wilf-equivalent permutation
matrices that are not equirestrictive.

Lastly, let us mention that all the discussion of this section can
be carried out with almost no changes to the case where the matrix
we want to avoid in the filling can have arbitrary nonnegative
entries; this corresponds to avoiding split graphs with multiple
edges. The interested reader will have no problems in filling in
the details.


\section{Proof of Theorem~\ref{thm:avoidk}}\label{sec:proof}


This section is devoted to the proof of Theorem~\ref{thm:avoidk}.
We show that we can adapt to our setting the proof
of~\cite{filling_boards}, which is for shape-Wilf-equivalence,
that is, row and column sums equal to $1$; we include the details
for the sake of completeness. (Actually, \cite{filling_boards}
contains two proofs of the analogous of our
Theorem~\ref{thm:avoidk} for shape-Wilf-equivalence; the proof we
adapt is the first one.) This bijection has been further studied
in~\cite{bms}. Here we show that it extends, in a quite
straightforward way, to arbitrary fillings. This gives a result
stronger than Theorem~\ref{thm:avoidk}, the consequences of which
in graph theoretic terms have already been pointed out at the end
of the previous section. Let us also mention that
Theorem~\ref{thm:avoidk} can also be proved using the techniques
of~\cite{krattenthaler}.

From now on $\bar{T}$ denotes a diagram with prescribed row and
column sums.
 When we say that a cell is
above (or below, to the right, to the left) of another cell we
always mean strictly. If we say that a cell is weakly above
(below, etc.) we mean not above (not below, etc.)

If $A$ and $B$ are two matrices, by $[A|B]$ we mean the matrix having
$A$ and $B$ as blocks, that is,
$$\left( \begin{array}{cc} A & 0 \\ 0 & B\end{array}\right) .$$

\begin{prop}\label{prop:avoidextension}
Let $M$ and $N$ be a pair of equirestrictive matrices and let $A$
be any matrix. Then the matrices $[M|A]$ and $[N|A]$ are also
equirestrictive.
\end{prop}

\begin{proof}
Let $L$ be a filling of the diagram $\bar{T}$ that avoids $[M|A]$.
Let $T'$ be the set of cells $(i,j)$ of $T$ such that the cells to
the right and below $(i,j)$ contain the matrix $A$. $T'$ is a
diagram, since if $(i,j)$ is in $T'$ all the cells weakly above
and weakly to the left of it are also in $T'$.
 Now set the row and column sums of $T'$ according to
 the restriction of $L$ to $T'$, call it $L'$, giving a diagram $\bar{T'}$. Now
$L'$ is a filling of $\bar{T'}$ that avoids $M$, so by assumption there is a
bijection between such fillings and the ones that avoid $N$. Change the
entries of $L$ corresponding to $T'$ to obtain a filling  of $\bar{T}$ that
avoids $[N|A]$.

The bijection in the other direction goes just in the same way.
\end{proof}

Let $F_t$ be the  matrix $[J_{t-1}|I_1]$. The proof of the following
proposition takes the rest of this section.

\begin{prop}\label{prop:jt_ft}
For all $t$, $F_t$ and $J_t$ are equirestrictive.
\end{prop}

We get as a corollary a stronger version of Theorem~\ref{thm:avoidk}.

\begin{cor}\label{cor:it_jt}
For all $t$, $[I_t|A]$ and $[J_t|A]$ are equirestrictive.
\end{cor}

\begin{proof}
By Proposition~\ref{prop:avoidextension} it is enough to show that
$I_t$ and $J_t$ are equirestrictive. The proof is by induction on
$t$; clearly $I_1$ and $J_1$ are equirestrictive. By
 Proposition~\ref{prop:jt_ft}, it is enough to show that $I_{t}$ and
$F_{t}$ are equirestrictive, and this follows by the induction
hypothesis combined with Proposition~\ref{prop:avoidextension}.
\end{proof}

A sketch of the proof of Proposition~\ref{prop:jt_ft} is as follows.
We first define two maps between fillings that
transform occurrences of $F_t$ into occurrences of $J_t$, and conversely, and
use them to define to algorithms that transform a filling avoiding
$F_t$ into a filling avoiding $J_t$, and conversely. The fact that these
two algorithms are inverses of each other follows from a series of lemmas.

For any filling $L$,
given two occurrences $G_1$ and $G_2$ of $J_t$ in $L$, we say that $G_1$
\emph{precedes} $G_2$ if the first entry in which they differ,
from left to right,
is either higher in $G_1$ or it is at the same height and the one in $G_1$ is
 to the left. So two occurrences are either equal or comparable.

The order for the occurrences of $F_t$ goes the other way around, i.e., we
look at the first entry in which they differ, from right to left, and the
lower entries have preference, and if they are at the same height, the one
more to the right goes first.

Let $L$ be a filling with the first occurrence of $J_t$ in rows
$r_1,\ldots,r_t$ and columns $c_1,\ldots,c_t$. Let $\phi(L)$ be
the result of substracting $1$ from each cell $(r_s,c_s)$, $1\leq s\leq t$ and
adding $1$ to each cell $(r_s,c_{s-1})$, $2\leq s\leq t$
and to cell $(r_1,c_t)$.
Since row and column sums have not been altered,   $\phi(L)$ is a
filling of $\bar{T}$. So we have changed an occurrence of $J_t$ to
an occurrence of $F_t$. Define $\psi$ as the inverse procedure,
that is, $\psi$ takes a filling of the diagram, looks for the
first occurrence of $F_t$, and replaces it by an occurrence of
$J_t$.

We define the algorithms $A1$ and $A2$ in the following way.
Algorithm $A1$ starts with a filling avoiding $F_t$ and applies
$\phi$ successively until there is no occurrence of $J_t$. The
result (provided the algorithm finishes) is a filling that avoids
$J_t$. Similarly, algorithm $A2$ starts with a filling avoiding
$J_t$ and applies $\psi$ until there are no occurrences of $F_t$
left. We claim that $A1$ and $A2$ are inverse of each other. We
prove this through a series of analogous lemmas.  It is enough to
prove the following claims.

\begin{itemize}
\item That both algorithms end. (Lemmas~\ref{lem:nojtabove}
and~\ref{lem:noftbelow2}.)
\item That $\psi(\phi^n(L))=\phi^{n-1}(L)$ for all $n$.
(Lemma~\ref{lem:psipicksb}.)
\item That $\phi(\psi^n(L))=\psi^{n-1}(L)$ for all $n$.
(Lemma~\ref{lem:phipicksa}.)
\end{itemize}

In order to prove these
claims, we need to investigate some properties of the maps $\phi$
and $\psi$. We start by studying the map $\phi$.

Let us first introduce some notation. Let $L$ be a filling of the
diagram and let $a_1,\ldots,a_t$ be the cells of the first $J_t$
in $L$, listed from left to right; say they are
$(r_1,c_1),\ldots,(r_t,c_t)$. So in each cell $a_i$ there is a
positive integer, possibly greater than one. Let $b_1,\ldots,b_t$
be the cells $(r_2,c_1),(r_3,c_2),\ldots,(r_t,c_{t-1})$ and
$(r_1,c_t)$; hence, $b_1,\ldots,b_t$ are the cells corresponding
to the occurrence of $F_t$ that is created after applying $\phi$
to $L$. So cell $b_i$ is in the same row as $a_{i+1}$ and in the
same column as $a_i$, for $i$ with $1\leq i \leq t-1$.

\begin{figure}[ht]
\begin{center}
\includegraphics[width=9cm]{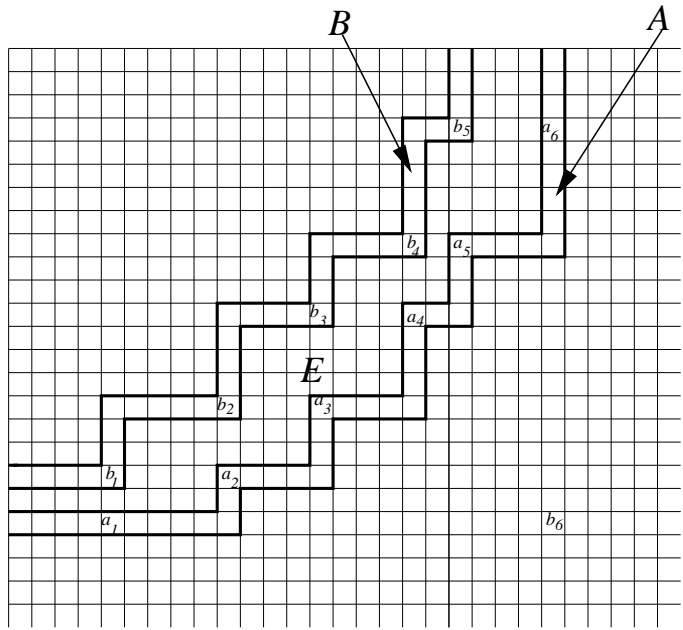}
\caption{The regions $A, B,$ and $E$}\label{fig:abe}
\end{center}
\end{figure}

 Consider now the following two paths of cells
determined by $a_1,\ldots,a_t$ and $b_1,\ldots,b_t$ (see
Figure~\ref{fig:abe}). The path $A$ starts at the leftmost cell in
the row of $a_1$, continues to the right until it reaches the
column of $a_2$, then takes this column up until it hits cell
$a_2$, then turns right until reaching the column of $a_3$, goes
up until $a_3$, then turns right again, and so on, until it
reaches cell $a_t$, at which point continues up until the top of
the diagram. The path $B$ is defined in a similar manner. It
starts at the leftmost cell of the row of cell $b_1$, and goes
right until it hits $b_1$. Then it turns up until the row of
$b_2$, where it turns and continues to the right until hitting
$b_2$. Then it goes up until the row of $b_3$, and then turns to
the right until $b_3$, and so on, until reaching $b_{t-1}$, at
which point it goes up until reaching the top of the diagram.
Since
 $a_1,\ldots,a_t$ are the first occurrence of $J_t$,
 the cells that are both to the right
of $B$ and to the left of $A$ are empty, or, in other words, this
region of the diagram avoids $J_1$. We denote this region  by $E$.
The choice of the first $J_t$ also imposes some other less trivial
bounds on the longest $J_i$'s that can be found in some other
areas determined by $E$. Note that in the next lemma the area left
of $E$ includes the path $B$.

\begin{lemma}\label{lem:leftofE}
With the above notation, the following hold for any filling $L$
and for the corresponding $\phi(L)$.
\begin{itemize}
\item[(i)] For all $i$ with $1\leq i \leq t-1$, there is no $J_i$
below $b_i$ and to the left of $E$.

\item[(ii)] For all $i$ with $1\leq i \leq t-1$, there is no
$J_{t-i}$ above and to the right of $b_i$ and to the left of $E$.

\item[(iii)] For all $i,j$ with $1\leq i<j\leq t-1$, the rectangle determined
by $b_i$ and $b_j$ contains no $J_{j-i}$ to the left of $E$; that is,
there is no $J_{j-i}$
 below $b_j$, above $b_i$, to the right of $b_i$, and to the left of $E$.
\end{itemize}
\end{lemma}

\begin{proof}
The arguments below apply to both $L$ and $\phi(L)$ since they do
not use the entries in cells $b_i$.

\begin{itemize}
\item[(i)] Assume there was such a $J_i$. Then this $J_i$ together
with $a_{i+1},\ldots,a_t$ would form a $J_t$ contradicting the
choice of $a_1$.

\item[(ii)] Suppose there was such a $J_{t-i}$. Then
$a_1,\ldots,a_i$ followed by this $J_{t-i}$ form a $J_t$ that
contradicts the choice of $a_{i+1}$.

\item[(iii)] Again, if there was such a $J_{j-i}$, combined with
$a_1,\ldots,a_{i-1}$ and $a_{j+1},\ldots,a_t$, it would create a
$J_t$ contradicting the choice of $a_i$.
\end{itemize}
\end{proof}

\begin{lemma}\label{lem:nojtabove}
There is no $J_t$ in $\phi(L)$ in the rows above $a_1$
\end{lemma}

\begin{proof}
We argue by contradiction. Let $G$ be an occurrence of
  $J_t$ in $\phi(L)$. Since
$\phi$ picked $a_1$ as the topmost cell being the left-bottom cell
 of a
$J_t$, $G$ must use at least one of the cells
$b_1,\ldots,b_{t-1}$. The idea is to substitute these cells $b_i$, and
possibly others, by some of the cells $a_i$, to find an occurrence of $J_t$
in $L$ in the rows above $a_1$,  hence contradicting the choice of $a_1$.

Now for each cell $b_{i}$ which belong to $G$, find the largest integer
$j$ such
that all cells of $G$ above $b_{i}$ and weakly below $b_{j-1}$ lie left
of $E$. In this way it is possible to find two sequences $i_1,\ldots,i_s$ and
 $j_1,\ldots, j_s$ with the following properties:
\begin{itemize}
\item $i_k<j_k$, $1\leq i_{k-1}<i_k$, and $j_{k-1}<j_k\leq t$ for all $k$;
\item $b_{i_k}$ is in $G$;
\item if $b_l$ is in $G$, then $i_k\leq l \leq j_k-1$ for some $k$;
\item all cells of $G$ above $b_{i_k}$ and weakly below $b_{j_k-1}$ are to
the left of $E$, and $j_k$ is the largest integer with this property.
\end{itemize}
Now we show that we can replace the cells of $G$ that fall left of $E$ and are
contained in the rectangles
determined by $b_{i_k}$ and $b_{j_k}$ by some of the $a_i$, giving an
instance of $J_t$ contained in $L$ and above $a_1$. We need to distinguish two cases,
according to whether $j_s=t$ or not. Assume first that $j_s\neq t$. For each
$k$, consider
the rectangles determined by $b_{i_k}$ and $b_{j_k}$. By
Lemma~\ref{lem:leftofE}.(iii), there are at most $j_k-i_k-1$ elements of
$G$ in this rectangle and to the left of $E$. Replace these cells, together
with $b_{i_k}$, by a (possibly proper) subset of $a_{i_k+1},\ldots,a_{j_k}$.
After doing this for each $k$, we still have an occurrence of $J_t$ starting
above $a_1$, but now it is contained in the original filling $L$, contradicting
the hypothesis. Now assume that $j_s=t$. For $k<s$, do the same substitutions
as in the previous case; for $k=s$, we have by Lemma~\ref{lem:leftofE}.(ii)
that there are at most $t-i_k-1$ cells of $G$ left of $E$ and above $b_{i_k}$.
Replace these cells and $b_{i_k}$ by a subset of $a_{i_k+1},\ldots,a_t$.
Again we obtain an occurrence of $J_t$ in $L$ that starts above $a_1$,
a contradiction.

 \end{proof}

This lemma alone shows that algorithm $A1$ terminates. Indeed, after one
application of $\phi$ all the cells in the row of $a_1$ and to the left
of $a_1$ are empty (because of the choice of $a_1$),
 and the cell $a_1$ has decreased its value by one. So
the leftmost cell of the first occurrence of $J_t$ in $\phi(L)$ is either
$a_1$, or it is to the right of $a_1$, or it is below $a_1$.
 But since the value in
cell $a_1$ decreases and cells to the left of $a_1$ stay empty,
 eventually there will be no occurrence of $J_t$
whose leftmost cell is $a_1$. So the selection of $J_t$'s goes from top to
bottom and from left to right, so for some $n$ the filling $\phi^n(L)$ is
 free of $J_t$'s.

It is not the case that if we apply $\phi$ to  an arbitrary filling $L$ of
$\bar{T}$ we have that $\psi(\phi(L))=L$. But algorithm $A1$ starts with a
filling that avoids $F_t$ and the successive applications of $\phi$ create
occurrences of $F_t$ from top to bottom and from left to right. We need
to show that in this situation
after each application of $\phi$, the first occurrence
of $F_t$ is precisely the one created by $\phi$. The next lemmas are
devoted to proving this.

\begin{lemma}\label{lem:noftbelow}
If $L$ contains no $F_t$ with at least one square below
$a_1$, then $\phi(L)$ contains no such $F_t$.
\end{lemma}

\begin{proof}
The proof is similar to the one of the previous lemma. Let $G$ be an
occurrence of $F_t$ in $\phi(L)$ with at least one cell below $a_1$. Since
$L$ had no such occurrence, $G$ contains at least one of the cells $b_i$.
 The bottom-right
cell of $G$ is below $a_1$, and it cannot be to the right of $a_{t-1}$,
otherwise this cell together with $a_1,\ldots, a_{t-1}$ would form an $F_t$
in $L$. By an argument similar to the one in the previous lemma, we change
all cells $b_i$ of $G$, and possibly others, to some of the cells $a_i$,
so that at the end we have an occurrence of $J_{t-1}$ that together with
the bottom-right cell of $G$ gives an occurrence of $F_t$ that contradicts
the hypothesis.

For each $b_i$ that is in $G$, look for the smallest $j$ such that all
cells in $G$ that are left of $b_i$ and weakly to the right of $b_{j+1}$ are
left of $E$. By doing this we find integers
$i_1,\ldots, i_s$ and $j_1,\ldots, j_s$ with the following
properties:
\begin{itemize}
\item $i_k> j_k$, $t-1\geq i_{k-1}>i_k$, and $j_{k-1}>j_k\geq 0$ for all $k$;
\item $b_{i_k}$ is in $G$ for all $k$ with $1\leq k \leq s$;
\item $j_k$ is the smallest integer such that all cells of $G$ that are left of
 $b_{i_k}$ and weakly to the right of $b_{j_k+1}$ are to the left of $E$;
\item if $b_l$ is in $G$, then $j_k+1 \leq l\leq i_k$ for some $k$.
\end{itemize}

We have to distinguish whether
${j_s}=0$ or not. Assume first $j_s\neq 0$. Since by
Lemma~\ref{lem:leftofE}.(iii) there are at most $i_k-j_k-1$ cells of $G$
in the rectangle determined by $b_{j_k}$ and $b_{i_k}$, these cells,
together with $b_{i_k}$,
can be replaced by a (possibly proper) subset of
$a_{j_k+1},\ldots,a_{i_k}$. By doing this for all $k$, we have an occurrence
of $J_{t-1}$ in $L$ that together with the right-bottom cell of $G$
contradicts the hypothesis. If $j_s=0$, then we do the same substitutions
for all $k\neq s$; for $k=s$, we have by Lemma~\ref{lem:leftofE}.(i) that
there are at most $i_s-1$ cells of $G$ left of $E$ and below $b_{i_s}$, so
we can substitute those and $b_{i_s}$ by $a_1,\ldots,a_{i_s}$. After these
substitutions, the result is again an occurrence of $F_t$ in $L$ that contains
a cell below $a_1$, contradicting the hypothesis.
\end{proof}

The following is easy but we state it for the sake of completeness.

\begin{lemma}\label{lem:noftright}
If $L$ contains no $F_t$ with a cell to the right of $a_t$ and below $a_2$, then $\phi(L)$
contains no such $F_t$.
\end{lemma}

\begin{proof}
Again we argue by contradiction. Suppose $G$ is an $F_t$ in $\phi(L)$
that contains a cell to the right of $a_t$ and below $a_2$. This cell
together with $a_2,\ldots,a_t$ gives an occurrence of $F_t$ in $L$ that
contradicts to the assumption.
\end{proof}

\begin{lemma}\label{lem:ftfromphi}
For each $k$ with $1\leq k\leq t-1$, there is no $J_k$ in $\phi(L)$
above $a_1$ and to the
left of and below $a_{k+1}$.
\end{lemma}

\begin{proof}
Let $G$ be an occurrence of such a $J_k$. If $G$ contains
none of $b_1,\ldots, b_{k-1}$, then $G$ followed by $a_{k+1},\ldots,a_t$
forms a $J_t$ in $L$ that is above $a_1$, and this contradicts the
choice of $a_1$. Hence, $G$ uses some $b_i$ for $1\leq i \leq k-1$.
 By an argument analogous to that of the proof of
Lemma~\ref{lem:nojtabove}, we can substitute the cells $b_i$ that are in $G$
and possibly others by some
$a_i$'s so that we get an occurrence of $J_{k}$ in $L$ that is below
$a_{k+1}$ and above $a_1$. This followed by $a_{k+1},\ldots,a_t$,
gives an $J_t$ in $L$ that contradicts the choice of $a_1$.
\end{proof}

The following lemma is just a combination of the previous and
induction; it implies that the inverse of algorithm $A1$ is $A2$.

\begin{lemma}\label{lem:psipicksb}
\begin{itemize}
\item[(i)] If $L$ does not contain any occurrence of $F_t$ below $a_1$, then the
first
occurrence of $F_t$ in $\phi(L)$ is $b_1,\ldots,b_t$.
\item[(ii)] If $L$ is a filling that avoids $F_t$, then
$\psi(\phi^n(L))=\phi^{n-1}(L)$.
\end{itemize}
\end{lemma}

\begin{proof}
For the first statement, let $f_1,\ldots,f_t$ be the first occurrence of $F_t$ in $\phi(L)$, with
the elements ordered from left to right. Recall that $b_1,\ldots,b_t$ is
an occurrence of $F_t$ in $\phi(L)$; we need to show that $f_i=b_i$
for all $i$.  By Lemma~\ref{lem:noftbelow},
$f_t$ is in the same row as $b_t$. By Lemma~\ref{lem:noftright}, $f_t$
cannot be to the right of $a_t$, hence $f_t=b_t$. Now use induction on
$t-i$. Suppose we know  $f_{i+1}=b_{i+1},\ldots,f_t=b_t$. It is enough now to
show that $f_i$ lies in the same row as $b_i$, since all
the cells to the right of $b_i$ but left of $b_{i+1}$ lie in $E$, which
we know contains only empty cells. But now Lemma~\ref{lem:ftfromphi} guarantees
that there is no $J_i$ below $b_i$, to the left of $b_{i+1}$, and
above $b_t$, as required.

For the second statement, it follows by Lemma~\ref{lem:noftbelow} and
induction on $n$ that the filling $\phi^{n}(L)$ contains no $F_t$ whose
lowest cell is below the lowest cell of the first occurrence of $J_t$. Hence
the previous statement applied to $\phi^n(L)$ gives immediately
that $\psi(\phi^n(L))=\phi^{n-1}(L)$.
\end{proof}

So the inverse of algorithm $A1$ is $A2$. Now
we only need to prove the converse. The proof follows exactly the same
steps and we content ourselves by stating and proving the corresponding
lemmas. Actually in this case some proofs are slightly simpler.

We keep the notation as above. Let $L$ be now a filling of $\bar{T}$ and let
$b_1,\ldots,b_t$ be the first occurrence of $F_t$ and let $a_1,\ldots,
a_t$ be the  occurrence of $J_t$ in $\psi(L)$ created after applying $\psi$ to
$L$. Consider again
the region $E$ as defined above. By the choice of $b_1,\ldots,b_t$ as the
first occurrence of $F_t$ in $L$, all the
cells of $E$ are again empty.

\begin{lemma}\label{lem:rightofE}
For all $i,j$ with $1\leq i<j\leq t$, the rectangle determined
by $a_i$ and $a_j$ contains no $J_{j-i}$ to the right of $E$ in either
$L$ or $\psi(L)$; that is,
there is no $J_{j-i}$
 below $a_j$, above $a_i$, to the left of $a_j$, and to the right of $E$.
\end{lemma}

\begin{proof}
Suppose there was such a $J_{j-i}$. Then $b_1,\ldots,b_{i-1}$, followed
by this $J_{j-i}$ and then followed by $b_j,\ldots,b_{t}$ gives an occurrence
of $F_t$ in $L$ that contradicts the choice of $b_{j-1}$.
\end{proof}

\begin{lemma}\label{lem:noftbelow2}
There is no $F_t$ in $\psi(L)$ with at least one cell in a row below $a_1$.
\end{lemma}

\begin{proof}
Suppose there is such an $F_t$.
 Its right-bottom cell is below $a_1$ and
also weakly to the left of $b_{t-1}$, since otherwise $b_1,\ldots,b_{t-1}$
and this cell would form an $F_t$ contradicting the choice of $b_t$. Let $G$ be
this occurrence of $F_t$ except the right-bottom cell.
$G$ must contain some of the cells $a_1,\ldots,a_t$. As in the previous lemmas,
the idea is to substitute the $a_i$ in $G$ together with other cells
 by some of the $b_i$ so that we
obtain an occurrence of $F_t$ in $L$ contradicting the choice of $b_t$.
Find integers $i_1,\ldots,i_s$ and $j_1,\ldots,j_s$ with the
following properties:

\begin{itemize}
\item $ i_k< j_k$, $1\leq i_{k-1}<i_k$, and $j_{k-1}<j_k\leq t-1$ for all $k$;
\item $a_{i_k}$ is in $G$ for all $k$ with $1\leq k \leq s$;
\item $j_k$ is the largest integer such that all cells of $G$ that are
to the right of $a_{i_k}$ and weakly to the left of $a_{j_k-1}$
are to the right of $E$;
\item if $a_l$ is in $G$, then $i_k\leq l \leq j_k-1$ for some $k$.
\end{itemize}

Now, by Lemma~\ref{lem:rightofE}, there are at most $j_k-i_k-1$ elements
of $G$ in the rectangle determined by $a_{i_k}$ and $a_{j_k}$. Together
with $a_{i_k}$,
they account for at most $j_k-i_k$ elements of $G$; substitute them for
a subset of
$b_{i_k},\ldots,b_{j_k-1}$. Doing this for all $k$, we get an occurrence of
$F_t$ in $L$ that contains a cell below $a_1$, hence contradicting the
choice of $b_1,\ldots,b_t$ as the first $F_t$ in $L$.
\end{proof}

\begin{lemma}\label{lem:nojtabove2}
If $L$ contains no $J_t$ that is above $a_1$, then $\psi(L)$ contains
no such $J_t$.
\end{lemma}

\begin{proof}
Let $G$ be such a $J_t$; $G$ must contain some of the cells $a_i$.
Find integers $i_1,\ldots,i_s$ and $j_1,\ldots,j_s$ with the
following properties:

\begin{itemize}
\item $ i_k > j_k$, $t\geq i_{k-1}>i_k$, and $j_{k-1}>j_k\geq 1$ for all $k$;
\item $a_{i_k}$ is in $G$ for all $k$ with $1\leq k \leq s$;
\item $j_k$ is the smallest integer such that all cells of $G$ that are
below $a_{i_k}$ and weakly above $a_{j_k+1}$ are to the right of $E$;
\item if $a_l$ is in $G$, then $j_k+1\leq l\leq i_k$ for some $k$.
\end{itemize}

As in the proof of the previous lemma, it is possible to substitute the
elements of $G$ contained in the rectangles determined by $a_{i_k}$ and
$a_{j_k}$, plus the cell $a_{i_k}$,
by (a subset of) the elements $b_{j_k},\ldots,b_{i_k-1}$. These substitutions
give a $J_t$ in $L$ that is above $a_1$, contrary to the hypothesis.
\end{proof}

\begin{lemma}\label{lem:nojtleft}
If $L$ contains no $J_t$ with a cell to the left of $a_1$ and below $a_2$,
then neither does $\psi(L)$.
\end{lemma}

\begin{proof}
If this were the case, the leftmost cell of this $J_t$ together with
$b_1,\ldots,b_{t-1}$ would give a $J_t$ contradicting the hypotheses
\end{proof}

\begin{lemma}\label{lem:jtfrompsi}
If $L$ contains no $J_t$ above $a_1$, there is no $J_{t-r}$  in $\psi(L)$
above  $a_{r+1}$ such
that the lowest cell of this $J_{t-r}$ is weakly to the left of $a_{r+1}$.
\end{lemma}

\begin{proof}
Suppose $G$ is an occurrence of such a $J_{t-r}$. $G$ must contain some of the cells $a_{r+2},
\ldots,a_t$, otherwise $b_1,\ldots,b_r$ followed by $G$ would form a $J_t$
contradicting the hypothesis.
Find integers $i_1,\ldots,i_s$ and $j_1,\ldots,j_s$ with the
following properties:

\begin{itemize}
\item $ i_k > j_k$, $t\geq i_{k-1}>i_k$, and $j_{k-1}>j_k\geq r-1$ for all $k$;
\item $a_{i_k}$ is in $G$ for all $k$ with $1\leq k \leq s$;
\item $j_k$ is the smallest integer such that all cells of $G$ that are
below $a_{i_k}$ and weakly above  $a_{j_k+1}$ are to the right of $E$;
\item if $a_l$ is in $G$, then $j_k+1\leq l\leq i_k$ for some $k$.
\end{itemize}

As before, the rectangle determined by $a_{i_k}$ and $a_{j_k}$ contains
at most $i_k-j_k-1$ cells of $G$; these cells, together with $a_{i_k}$,
 can be replaced by a subset of
$b_{j_k},\ldots,b_{i_k-1}$. After all these substitutions we get an occurrence
of $J_{t-r}$ in $L$ that combined with $b_1,\ldots,b_r$ gives an occurrence
of $J_t$ in $L$ contradicting the hypothesis.
\end{proof}

\begin{lemma}\label{lem:phipicksa}
\begin{itemize}
\item[(i)] If $L$ does not contain any occurrence of $J_t$ above $b_t$, then the
first
occurrence of $J_t$ in $\phi(L)$ is $a_1,\ldots,a_t$.
\item[(ii)] If $L$ is a filling that avoids $J_t$, then
$\phi(\psi^n(L))=\psi^{n-1}(L)$.
\end{itemize}
\end{lemma}

\begin{proof}
For the first statement, let $d_1,\ldots,d_t$ be the first occurrence of
$J_t$ in $\psi(L)$, with cells listed from left to right. We want to show
that $a_i=d_i$ for all $i$ with $1\leq i \leq t$. By
Lemma~\ref{lem:nojtabove2}, $d_1$ is in the same row as $a_1$, and by
Lemma~\ref{lem:nojtleft} it is weakly to the right of $a_1$, hence $d_1=a_1$.
Now we proceed by induction on $i$. Suppose $d_1=a_1,\ldots,d_i=a_i$.
By Lemma~\ref{lem:jtfrompsi} we have that the only $J_{t-i}$ in $\psi(L)$
 that is weakly above and weakly to the left of $a_{i+1}$ is $a_{i+1},\ldots,
a_t$, hence $d_{i+1}=a_{i+1}$, as needed.

For the second statement, by induction and Lemma~\ref{lem:nojtabove2} we
get that $\psi^n(L)$ satisfies the hypothesis of part (i), hence it follows
that  $\phi(\psi^n(L))=\psi^{n-1}(L)$.
\end{proof}


\section{Concluding remarks}\label{sec:comments}


In his paper~\cite{krattenthaler}, Krattenthaler speaks of a
``bigger picture" that would englobe several recent results on
pattern avoiding fillings of diagrams. We believe that our
correspondence between graphs and fillings of diagrams also
belongs to this picture and that it may shed some light in the
understanding of it. We have shown that for each statement in
pattern avoiding fillings there is a statement about graphs
avoiding certain split graphs. So we can claim that in some sense
the resources available to attack either problem have doubled. An
example of this are the ``repeated" results in the literature
mentioned at the end of Section~\ref{sec:degree}.

For completeness, we mention here a result by Bousquet-M\'elou and
Steingr\'imsson~\cite{bms} that can be cast in terms of
$k$-noncrossing and $k$-nonnesting graphs. They restrict to
diagrams with self-conjugate shape and row and column sums are set
to $1$, and they only consider symmetric $0-1$ fillings (that is,
symmetric with respect to the main diagonal of the diagram). For
these fillings, they show that $I_t$ and $J_t$ are
equirestrictive. In terms of matchings, this says that for each
left-right degree sequence, the number of $k$-noncrossing
symmetric matchings is the same as the number of $k$-nonnesting
ones, where a matching on $[2n]$ is symmetric if it equals its
reflection through the vertical axis that goes between vertices
$n$ and $n+1$. Similar results for symmetric graphs can be deduced
from~\cite[Theorem 15]{krattenthaler}.

\comment{

It is also possible to go in the other direction, that is, to use
graph theoretic results to find shape-Wilf-equivalent matrices.
The following is Theorem~1 of~\cite{match_and_part}, written in
our notation.

\begin{thm}\label{thm:thm1}
For any left-right degree sequence $D$ having only terms of the
form $(1,0),(0,1),$ and $(1,1)$, the number of graphs $G$ on $D$
with $\cross(G)=i$ and $\nest(G)=j$ is the same same as the number
of such graphs with $\cross(G)=j$ and $\nest(G)=i$.
\end{thm}

If we sum for $i,j$ with $1\leq i <k$ and $1\leq j <l$, we deduce
that the number of $k$-noncrossing and $l$-nonnesting graphs on
$D$ is the same as the number of $l$-nonnesting and
$k$-noncrossing such graphs. (This is the same as~\cite[Corollary
2]{match_and_part}, but fixing the degrees.)
 To rephrase
this in terms of fillings of diagrams, extend the previous
notation by saying that two sets of permutation matrices
$M_1,\ldots,M_s$ and $N_1,\ldots,N_r$ are equirestrictive if for
each diagram $T$ with prescribed row and column sums, the number
of fillings avoiding each of $M_1,\ldots,M_s$ equals the number of
fillings avoiding each of $N_1,\ldots,N_r$. If the row and column
sums are equal to $1$, we say the two sets of permutation matrices
are shape-Wilf-equivalent. Then as a corollary of
Theorem~\ref{thm:thm1} we get the following. (This is also
mentioned in the final section of~\cite{krattenthaler}.)

\begin{thm}\label{thm:i_k_j_l}
For all $k,l$, $\{I_k,J_l\}$ and $\{I_l,J_k\}$ are
shape-Wilf-equivalent.
\end{thm}

}

 Let us finish by going  back to our initial motivation of
studying  $k$-noncrossing and $k$-nonnesting graphs.  Even if our
main question has been answered positively, it is fair to say that
it has not been solved in the most satisfactory way; ideally we
would like to find a bijective proof in graph theoretic terms.
Note that due to its roundabout character, our proof of
Theorem~\ref{thm:avoidk} does not give a clear bijection, neither
in terms of graphs nor of fillings. Also the proofs of
Corollaries~\ref{cor:kcrossknestmulti}, \ref{cor:simple},
and~\ref{cor:simplemax} do not provide bijections in graph
theoretic terms. A bijective proof of
Theorem~\ref{thm:kcrossknest} for $k=2$ has recently been found by
Jel\'inek, Klazar, and de Mier~\cite{notes_jkm}.

Other interesting questions related to $k$-crossings and
$k$-nestings of graphs include, as mentioned before, to determine
whether the pairs $(\cross(G), \nest(G))$ are symmetrically
distributed among all graphs. This is already known for matchings
and partition graphs~\cite{match_and_part}. One would also hope
for a wide generalization of Theorem~\ref{thm:kcrossknest} stating
that the number of graphs with $r$ $k$-crossings and $s$
$k$-nestings equals the number of graphs with $s$ $k$-crossings
and $r$ $k$-nestings. Again, the case $k=2$ is known  for
matchings~\cite{klazar_match_trees} and partition
graphs~\cite{distribution}. Unfortuntely,  for $k=3$ this is not
true even for matchings; for instance, Marc Noy~\cite{marc}
 checked that there are more matchings with six edges and only one $3$-crossing than with only
one $3$-nesting.

\comment{ We do not get symmetry of the distribution of the number
of occurrences of $I_t$ and $J_t$ in fillings of a given
$\bar{T}$. An easy example is considering the diagram of shape
$(2,2)$ with column sums $(2,1)$ and row sums $(1,2)$. There are
only two fillings, by rows $(1,0|1,1)$ and $(0,1|2,0)$, the first
having one $I_2$ and the second having two $J_2$. If we only look
at the cells, and not the values, they do have the same number of
occurrences actually. hmm... is there a counterexample to
symmetric distribution if we count each occurrence by the minimum
value of the cells it occupies? actually, this makes sense in
terms of the algorithm... interesting...}


\section*{Acknowledgements}


I am very grateful to Sergi Elizalde, V\'it Jel\'inek, Martin
Klazar, Martin Loebl, and Marc Noy for many fruitful and
stimulating discussions and for pointing out several useful
references.


\end{document}